\def\e{{\epsilon}}

\def\w{{\omega}}
\def\W{{\Omega}}

\def\cA{{\cal A}}

\def\cC{{\cal C}}
\def\cD{{\cal D}}
\def\cF{{\cal F}}
\def\cG{{\cal G}}

\def\cM{{\cal M}}
\def\cN{{\cal N}}
\def\cP{{\cal P}}

\def\bE{{\rm\bf E}}
\def\bN{{\rm\bf N}}

\def\t{{t-1}}

\def\inf{{\infty}}
\def\pf{{\diamond}}
\def\pr{{^\prime}}
\def\sh{{\partial}}

\def\lt{{\ \le\ }}

\def\pebb{\atopwithdelims<>}
\def\rar{{\ \rightarrow\ }}
\def\sub{{\ \subseteq\ }}
\def\sup{{\ \supseteq\ }}

\def\lf{{\lfloor}}
\def\rf{{\rfloor}}

\def\lgn{{\lg n}}

\def\sln{{\sqrt\lgn}}
\def\stln{{\sqrt{2\lgn}}}
\def\sn{{n^{1/2}}}
\def\sq{{\sqrt 2}}

\documentclass[12pt]{article}
\newtheorem{theorem}{Theorem}[section]

\newtheorem{lemma}[theorem]{Lemma}

\usepackage{amssymb}
\usepackage{latexsym,psfig,multicol}
\begin{document}

%######################################################################
%######################################################################
%
%       TITLE PAGE: 
%
\title{Thresholds for families of multisets,\\ with an application to graph pebbling}
\author{
Airat Bekmetjev\thanks{
Department of Mathematics,
Arizona State University,
Tempe, Arizona USA 85287-1804
(email: bekmetje@mathpost.la.asu.edu)}~\\
Graham Brightwell\thanks{
Department of Mathematics,
London School of Economics,
London, England WC2A 2AE
(email: g.r.brightwell@lse.ac.uk)}~\\
Andrzej Czygrinow\thanks{
Department of Mathematics,
Arizona State University,
Tempe, Arizona USA 85287-1804
(email: andrzej@math.la.asu.edu)}~\\
and\\
Glenn Hurlbert\thanks{
Department of Mathematics,
Arizona State University,
Tempe, Arizona USA 85287-1804
(email: hurlbert@math.la.asu.edu)}~\\
}
\maketitle
\newpage

%######################################################################
%######################################################################
%
%       ABSTRACT:
%
\begin{abstract}
\noindent
In this paper we prove two multiset analogs of classical results. 
We prove a multiset analog of Lov\' asz's version of the Kruskal-Katona Theorem 
and an analog of the Bollob\'as-Thomason threshold result. As a corollary we 
obtain the existence of pebbling thresholds for
arbitrary graph sequences. In addition, we improve both the lower and upper 
bounds for the `random pebbling' threshold of the sequence of paths.
\vspace{0.2 in}

\vspace{0.2 in}

\vspace{0.8 in}
\end{abstract}

\newpage

%######################################################################
%######################################################################
%
%       BEGINNING OF PAPER:
%
\section{Introduction}\label{Intro}
Throughout this paper $G$ will denote a simple connected graph, and
$n=n(G)$ will denote the number of its vertices.
The vertex set of $G$ will be the set $[n]=\{1,2,\ldots,n\}$.

%######################################################################
%       
%       SUBSECTION:
%       
\subsection{Pebbling}\label{Pebbling}
Suppose $t$ pebbles are distributed onto the vertices of a graph $G$.
A pebbling step $[v,v\pr]$ consists of removing two pebbles from one vertex $v$ 
and then placing one pebble at an adjacent vertex $v\pr$.
We say a pebble can be {\it moved} to a vertex $z$, the {\it root} vertex,
if we can repeatedly apply pebbling steps so that in the resulting
distribution $z$ has at least one pebble.

For a graph $G$, a distribution $D$ of pebbles onto the vertices of $G$, and a 
`root' vertex $z$, we say that $D$ is {\it $z$-solvable} if it is possible to
move a pebble to $z$; otherwise, $D$ is $z$-{\it unsolvable}.
Also $D$ is {\it solvable} if it is $z$-solvable for all $z$, and
{\it unsolvable} otherwise.

The {\it pebbling number} $pn(G)$ is the smallest integer $t$ such that 
{\it all} distributions of $t$ pebbles to the vertices of $G$ are solvable. 
In this paper we are concerned instead with the minimum $t$ such that 
{\it almost all} distributions of $t$ pebbles to the vertices of $G$ are 
solvable. 
The interested reader is encouraged to read \cite{Hur} for the history of
and the many results on graph pebbling.

%######################################################################
%       
%       SUBSECTION:
%       
\subsection{Random Distributions}\label{RandomDists}

In this paper we are interested in the probabilistic pebbling model, in which 
the pebbling distribution is selected uniformly at random from the set of all 
distributions with a prescribed number $t$ of pebbles (we emphasize that the 
pebbles are unlabeled and the vertices are labeled).  This is certainly not the 
only possible random pebbling model; for instance one could consider the model
obtained by placing each of the pebbles uniformly at random on a vertex of $G$
-- this model will exhibit very different behavior.  

We will study thresholds for the number $t$ of pebbles so that if
$t$ is essentially larger than the threshold, then a random
distribution is almost surely solvable, and if $t$ is essentially
smaller than the threshold, then a random distribution is almost surely 
unsolvable.
Formally, our notion of a pebbling threshold is defined as follows.
Let $\bN$ denote the set of nonnegative integers, and let
$D_n:[n]\rar \bN$ denote a distribution of pebbles on $n$ vertices.
For a particular function $t = t(n)$, we consider the probability space
$\cD_{n,t}$ of all distributions $D_n$ of {\it size} $t$, i.e.\ with
$t=\sum_{i\in[n]}D_n(i)$ pebbles, with each such distribution being equally 
likely.
Given a graph sequence $\cG=(G_1,\ldots,G_n,\ldots)$, where $G_n$ has vertex 
set $[n]$, denote by $P_\cG(n,t)$ the probability that an element of 
$\cD_{n,t}$ chosen uniformly at random is $G_n$-solvable.
We call a function $t=t(n)$ a {\it threshold} for $\cG$, and write 
$t\in th(\cG)$, if the following two statements hold for every sequence 
$\w=\w(n)$ tending to infinity (we write $\w\gg 1$): 
({\it i})~$P_\cG(n,t\w)\rar 1$ as $n\rar\inf$, and
({\it ii})~$P_\cG(n,t/\w)\rar 0$ as $n\rar\inf$.  
(Here and elsewhere, if $t\w$ and
$t/\w$ are not integers, they should be interpreted as taking on the nearest 
integer value.)  
Of course, the definition mimics the important threshold concept in random
graph theory. Unlike the situation in random graphs, however, it did not seem 
obvious that even ``natural'' families of graphs have pebbling thresholds, 
although the existence of a threshold for any graph sequence is conjectured 
in \cite{Hur}. 
The random pebbling model is also studied in \cite{CzyEHK}, where the following result 
is proved.

\begin{theorem}\label{gen}
For any $\e > 0$ and any graph sequence $\cG=(G_1,\ldots,G_n,\ldots)$, where 
$V(G_n)=[n]$ for each $n$, $th(\cG)\sub \Omega(\sn)\, \cap\ o(n^{1+\e})$.
\end{theorem}

In other words, {\it if} there is a threshold function $t(n)$ for a graph 
sequence, then $t(n)$ cannot be essentially smaller than $\sn$ (the threshold 
function for the family of complete graphs) or as large as $n^{1+\e}$.  
Our results in this paper serve to improve the upper bound.  

Better bounds were obtained in~\cite{CzyEHK} for graphs satisfying various 
diameter and connectivity conditions, and threshold functions were found 
for the classes of stars and wheels. 
Also in \cite{CzyEHK}, bounds were found for the classes of paths and cycles.

\begin{theorem}\label{paths}
\begin{description}
\item {\bf (a)}
For the sequence of paths $\cP=(P_1,P_2,\ldots,P_n,\ldots)$, 
we have $th(\cP)\sub \W(n)\cap o(n^{1+\e})$ for every $\e>0$.
\item {\bf (b)}
For the sequence of cycles $\cC=(C_1,C_2,\ldots,C_n,\ldots)$, 
we have $th(\cC)\sub$ $\W(n)\cap$ $o(n^{1+\e})$ for every $\e>0$.
\end{description}
\end{theorem}

Our purpose in this paper is firstly to show that every graph sequence does 
have a threshold, and secondly to prove much tighter bounds for the threshold 
functions of classes of paths and cycles.

\begin{theorem}
\label{SeqThresh}
Let $\cG=(G_1,\ldots,G_n,\ldots)$ be any graph sequence,
and define $t=t(n)=\min\{r \mid P_\cG(n,r)\ge 1/2\}$.
Then $t\in th(\cG)$.
\end{theorem}

\begin{theorem}\label{Thpaths}
For the sequence of paths $\cP=(P_1,\ldots,P_n,\ldots)$ we have
\begin{description}
\item{\bf (a)} $th(\cP)\sub\Omega(n2^{c\sln})$, where $c$ is any 
constant less than $1/\sq$;
\item{\bf (b)} $th(\cP)\sub O(n2^{2\sln})$.
\end{description}
\end{theorem}

Here and throughout, $\lg$ denotes the logarithm to base~2.  
Our proofs apply equally well to the sequence of cycles 
$\cC=(C_1,\ldots,C_n,\ldots)$.

In fact, Theorem \ref{SeqThresh} follows from a very general result about 
thresholds in random multiset models.  For a natural number $n$, let 
$\cM_n$ denote the partially ordered set ({\it poset}) of all multisets 
of $[n]$, ordered by inclusion.  For $\cF_n$ a subfamily of $\cM_n$, and $t$ a 
natural number, let $\cF_n(t)$ denote the family of $t$-element multisets in
$\cF_n$.  If $\cF_n\sub \cM_n$,   
the family $\cF_n$ is said to be {\it increasing} if $E\sup F\in\cF_n$
implies $E\in\cF_n$, and {\it decreasing} if $E\sub F\in\cF_n$ implies 
$E\in\cF_n$.

The size of $\cM_n(t)$ is ${n \pebb t} = {n+t-1 \choose t}$, and therefore 
$P_t(\cF_n(t)) \equiv |\cF_n(t)|/{n\pebb t}$ is the probability that a 
uniformly randomly chosen $t$-multiset of $[n]$ is in the family $\cF_n$.
We say that $t=t(n)$ is a {\it threshold} for a sequence 
$(\cF_1,\ldots,\cF_n,\ldots)$ of increasing families of multisets if, 
for any function $\w=\w(n)\gg 1$, we have $P_{t\w}(\cF_n(t\w))\rar 1$ and 
$P_{t/\w}(\cF_n(t/\w))\rar 0$.

\begin{theorem}\label{threxist}  
Let $\cF = (\cF_1,\dots,\cF_n,\ldots)$ be a sequence of increasing families, 
with $\cF_n \subseteq \cM_n$ for each $n$.  
Define $t = t(n) = \min \{ r \mid P_r(\cF_n(r)) \ge 1/2\}$.  
Then $t \in th(\cF)$. 
\end{theorem}

Theorem~\ref{threxist} is an analog of a result of Bollob\'as and 
Thomason~\cite{BolTho}, stating that any sequence of increasing families
of the {\it subset lattice} has a threshold function.    

One can view a particular distribution of $t$ pebbles to $n$ vertices as
a particular multiset of $t$ elements from the ground set $[n]$.
Thus the set of all solvable distributions for a particular graph on $n$
vertices corresponds to a family of multisets on $[n]$.
Since the set of all solvable distributions for a particular graph 
has the increasing property, Theorem \ref{SeqThresh} is an immediate 
consequence of Theorem \ref{threxist}.

In our proof of Theorem~\ref{threxist}, the crucial tool will be the following
1969 result of Clements and Lindstr\"om, extending an earlier result
of Macauley~\cite{Mac}, which is the multiset analog of the 
celebrated Kruskal-Katona Theorem~\cite{Fra,Kat,Kru} for the subset lattice. 

Given any subfamily $\cA\sub \cM_n(t)$, we define its {\it shadow}
$\sh\cA=\{C\in \cM_n(\t) \mid C\subset A\mbox{ for some } A \in\cA \}$, and set
$\sh^{i+1}\cA=\sh\sh^i\cA$.

For a multiset $A \in {\cal M}_n(t)$, and $i\in [n]$, let $A(i)$ denote
the multiplicity of $i$ in $A$.
The {\em colexicographic order} on ${\cal M}_n(t)$ is defined by setting
$A<B$ if $A\not=B$ and, for some $i \in [n]$, $A(i) < B(i)$ while
$A(j)=B(j)$ for $j>i$.

\begin{theorem} (Clements-Lindstr\"om)\label{CL}
Suppose that $\cF$ is a subset of $\cM_n(t)$, and that $\cG$ consists of the
first $|\cF|$ elements of $\cM_n(t)$ in colexicographic order.  Then, for any
$k \ge 1$,
$$
|\partial^k \cF| \ge |\partial^k \cG|. 
$$
\end{theorem} 

In other words, the size of the shadow (at any level) of a subset of $\cM_n(t)$
is minimized by taking an initial segment of the colexicographic order on
$\cM_n(t)$.

In 1979, Lov\'asz proved a version of the Kruskal-Katona theorem which was
used by Bollob\'as and Thomason~\cite{BolTho} to prove the existence of 
threshold functions.  
An analogous version of the Clements-Lindstr\"om theorem was
conjectured in \cite{Hur}.  We prove this conjecture in the next section.

For $x$ a non-negative real number, let ${x\pebb t}=(x)(x+1)\cdots (x+\t)/t!$.
(Note that this coincides with our earlier definition if $x$ is a natural 
number.)

\begin{theorem}\label{lov}
Suppose that $\cA\sub\cM_n(t)$ and define $x$ by $|\cA|={x\pebb t}$.
Then $|\sh\cA |\ge {x\pebb\t}$.
\end{theorem}

For the case in which $\cA = {m \pebb t}$ for $m$ a natural number, the first
$|\cA|$ elements of $\cM_n(t)$ in colexicographic order are the $t$-multisets
of $\{1,\dots, m\}$.  The shadow of the family consisting of these
multisets is the family of $(\t)$-multisets of $\{1,\dots,m\}$, of size 
${x \pebb \t}$, so Theorem~\ref{lov} is equivalent to Theorem~\ref{CL} in this
case.  For families of intermediate sizes, Theorem~\ref{lov} is a ``smoothed''
version of Theorem~\ref{CL}.  

The rest of this paper is organized as follows. In Section \ref{Sets} we prove
Theorem \ref{lov} (Section \ref{Lsec}) and Theorem \ref{threxist} 
(Section \ref{ExSec}).
Section \ref{PathThresh} is devoted to a proof of Theorem \ref{Thpaths}.

%######################################################################
%######################################################################
%       
%       NEW SECTION:
%       
\section{Set Theory}\label{Sets}
In this section, we prove the two results concerning the multiset lattice. 
First, verifying a conjecture from \cite{Hur}, we prove
our multiset analog of Lov\'asz's version of the Kruskal-Katona theorem. 
Second, we establish our multiset analog of the
Bollob\'as-Thomason threshold theorem \cite{BolTho}.  As an immediate corollary 
we obtain the existence of the pebbling threshold for any graph sequence.

\subsection{Multiset analog of Lov\'asz's Theorem}\label{Lsec}

Let $[n]\pebb t$ denote the family of $t$-element multisets of $[n]$,
and as before let ${n\pebb t} = {n+t-1 \choose t}$ denote its cardinality.
(For convenience in reading, the reader may enjoy using the terminology
``$n$ {\it pebble} $t$''.) 
Also as before, let $x\pebb t$ denote the polynomial evaluation of 
$n\pebb t$ for any real number $x$; that is, 
$${x\pebb t}=(x)(x+1)\cdots (x+\t)/t!.$$

Let $\cA\sub {[n]\pebb t}$, and for any $A\in\cA$ and $i\in [n]$, 
let $A(i)$ denote the multiplicity of $i$ in $A$.
We adopt the convention of writing an element of a multiset just once, with 
its multiplicity written as an exponent; for example $\{1^j\}$ denotes the
multiset of $j$ ones. Also, for any multisubset $I$ of $[n]$ we set 
$\cA-I=\{A-I\ |\ A\in\cA\}$.

In the proof of Theorem \ref{lov}, we make use of various different partitions 
of a family $\cA\sub {[n] \pebb t}$.  For $0\le j\le t$, define the sets 
$\cA_j=\{A\in\cA\ |\ A(1)=j\}$, forming a partition of $\cA$.
If $\cA ={[n]\pebb t}$, this partition gives rise to the relation
\medskip
\noindent
$${n\pebb t}\ =\ \sum_{j=0}^t\ {n-1\pebb j}\ ,$$
which has the polynomial equivalent
\begin{equation}\label{partn}
{x\pebb t}\ =\ \sum_{j=0}^t\ {x-1\pebb j}\ .
\end{equation}

For each $i\in [n]$, define the sets 
$\cA^i=\{A\in\cA\ |\ A(i)>0 {\rm\ and\ } A(j)=0 {\rm\ for\ all\ } j<i\}$, 
forming a second partition of $\cA$.
More important for us is that the sets $\cA_1,\ldots,\cA_t$ partition $\cA^1$.

A third partition of $\cA$ is given by $\cA =\cA^1\cup \cA_0$, 
which in the case that $\cA ={[n]\pebb t}$ gives rise to the relation
$${n\pebb t}\ =\ {n\pebb\t}\ +\ {n-1\pebb t}\ ,$$
having the polynomial equivalent
\begin{equation}\label{split}
{x\pebb t}\ =\ {x\pebb\t}\ +\ {x-1\pebb t}\ .
\end{equation}

Given a family $\cA \sub \cM_n(t)$, and indices $i,j$ with 
$1\le i<j\le n$, a {\it compression} of $\cA$ is obtained by taking
each member $A$ of $\cA$ such that $A(j)\ge 1$ and $A-\{j\}+\{i\} \notin \cA$,
and replacing it by $A-\{j\}+\{i\}$. 

A family $\cA$ is said to be {\it compressed} if, for all $1\le i<j\le n$, 
we have $A-\{j\}+\{i\}\in\cA$ whenever $A\in\cA$, i.e., $\cA$ is unchanged by
any compression.  Note that any family can be transformed into a compressed 
one by a sequence of compressions.  Note also that initial segments of the
colexicographic order are compressed families, but that these are not the only 
ones.   

If $\cA$ is compressed and $C\in\sh\cA$, then $C=A-\{i\}$ for some 
$A\in\cA$, $i\in [n]$.  Because $\cA$ is compressed, the set 
$E=A-\{i\}+\{1\}$ is in $\cA$.  Since $C=E-\{1\}$, we see that 
$\sh\cA\sub\cA^1-\{1\}$, which implies that
$$|\sh\cA |\lt |\cA^1|\ .$$
On the other hand, we know that $\cA^1-\{1\}\sub\sh\cA^1\sub\sh\cA$, and so
$$|\cA^1|\lt |\sh\cA |\ .$$
These two facts together imply the following lemma. 

\begin{lemma}\label{f1}
Let $\cA\sub\cM_n(t)$ and suppose that $\cA$ is a compressed family.
Then $|\sh\cA |=|\cA^1|$.
\hfill$\pf$
\end{lemma}

The following lemma is proved by Clements in \cite{Cle}.

\begin{lemma}\label{comp}
Suppose that $\cA\sub\cM_n(t)$ and let $q(\cA)$ be a compression 
of $\cA$.  Then $|\sh\cA|\ge |\sh q(\cA)|$.
\hfill$\pf$
\end{lemma}

With these tools we now can prove Theorem \ref{lov}.

\medskip
\noindent
{\it Proof of Theorem \ref{lov}.}
We use induction on $n$. Because of Lemma \ref{comp} we may assume that 
$\cA$ is compressed.
If $|\cA^1|\ge {x\pebb\t}$ then we are done because of Lemma \ref{f1}.
So we will assume that $|\cA^1|<{x\pebb\t}$ and argue to a contradiction.

\vspace{0.2 in}
\noindent
{\bf Claim\ }
For each $0\le j\le t$ we have $|\cA_j|\ge {x-1\pebb t-j}$.
\vspace{0.2 in}

The truth of this claim yields the following contradiction.
Because the families $\cA_1,\ldots,\cA_t$ partition the family $\cA^1$, we have
$$|\cA^1|\ 
	=\ \sum_{j=1}^t\ |\cA_j|\ 
	\ge\ \sum_{j=1}^t\ {x-1\pebb t-j}\ 
	=\ \sum_{k=0}^{\t}\ {x-1\pebb k}\
	=\ {x\pebb\t}\ ,$$
using equation (\ref{partn}).
Thus we need only to prove the claim.

Induction on $j$ shows that $|\sh^j\cA_0|\ge {x-1\pebb t-j}$.
Indeed, for $j=0$ we know from equation (\ref{split}) that 
$|\cA_0|>{x-1\pebb t}$.
Also, $\sh^j\cA_0$ is a family of $(t-j)$-multisets of the set $\{2,\ldots,n\}$
of size $n-1$, and so if $|\sh^j\cA_0|\ge {x-1\pebb t-j}$, then 
$|\sh^{j+1}\cA_0|\ge {x-1\pebb t-j-1}$ (by induction on $n$ for the theorem).

Now, if $C\in\sh^j\cA_0$ then $C=A-I$ for some $A\in\cA_0$ and some
submultiset $I$ of $A$ of size $j$.  Because $\cA$ is compressed, the set 
$E=A-I+\{1^j\}$ is in $\cA$, and in particular is in $\cA_j$.
Since $C=E-\{1^j\}$ we see that $\sh^j\cA_0\sub\cA_j-\{1^j\}$, which 
implies that 
$$|\cA_j|\ 
	=\ |\cA_j-\{1^j\}|\ 
	\ge\ |\sh^j\cA_0|\ 
	\ge\ {x-1\pebb t-j}\ .$$
This proves the claim, and the theorem follows.  
\hfill$\pf$
\medskip

%######################################################################
%######################################################################
%
%       SUBSECTION:
%
\subsection{Thresholds}\label{Thresholds}\label{ExSec}
The main result of this section is the analog of the Bollob\'as-Thomason 
threshold theorem.

\medskip
\noindent 
{\it Proof of Theorem \ref{threxist}.} We consider the two ``reference'' 
families: 
$$
\begin{array} {rcl}
\cM_n(r;b) &=& \{ A \in \cM_n(r) \mid A(n) < b\} \quad (1\le b \le r), \\
\cN_n(r;b) &=& \{ A \in \cM_n(r) \mid A(n-b+1) = \cdots = A(n) =0 \} \\
&& (1\le b \le n-1).
\end{array} 
$$
Note that each family $\cM_n(r;b)$ and $\cN_n(r;b)$ is an initial segment
of the colexicographic order on $\cM_n(r)$, and that, for any $k\ge 1$, 
$$
\partial^k \cM_n(r;b) = \cM_n(r-k;b)\ ; \qquad \partial^k \cN_n(r;b) 
= \cN_n(r-k;b)\ .
$$
We use these families, rather than general initial segments of the 
colexicographic order, since their sizes are a little easier to estimate. 
Our strategy is, for each $r$, to compare 
$\overline{\cF_n(r)}=\cM_n(r)-\cF_n(r)$ 
with an appropriate member of one of these reference families, 
at levels $r$ near $t$, somewhat above $t$, and somewhat below $t$.  

First, we need some estimates on the probabilities of the reference families, 
which we shall use repeatedly.  For any positive integers $n$, $r$ and $b$ 
with $b\le r$,  
$$
|\cM_n(r;b)| = |\cM_n(r)| - |\cM_n(r-b)| = {n \pebb r} - {n \pebb r-b}\ , 
$$
and so
$$
\begin{array} {rcl}
P_r(\cM_n(r;b)) &=& 1 - \frac{(n+r-b-1)!r!}{(n+r-1)!(r-b)!} \\
&=& 1 - \left( \frac{r}{n+r-1}\right) \left(\frac{r-1}{n+r-2}\right) \cdots 
\left(\frac{r-b+1}{n+r-b}\right).   
\end{array} 
$$
We derive the lower bounds
$$
\begin{array} {rcl}
P_r(\cM_n(r;b)) &\ge& 1 - \left(\frac{r}{n+r-1}\right)^b \\
&=& 1 - \left( 1 - \frac{n-1}{n+r-1}\right)^b \\
&\ge& 1 - \exp( -b(n-1)/(n+r-1))\ ,
\end{array}
$$
and the upper bounds
$$
\begin{array} {rcl}
P_r(\cM_n(r;b)) &\le& 1 - \left(\frac{r-b+1}{n+r-b}\right)^b \\
&=& 1 - \left( 1 + \frac{n-1}{r-b+1}\right)^{-b} \\
&\le& 1 - \exp( -b(n-1)/(r-b+1))\ . 
\end{array}
$$

Similarly we see that, for positive integers $n$, $r$ and $b$ with $b\le n-1$, 
$|\cN_n(r;b)| = |\cM_{n-b}(r)| = {n-b \pebb r}$, and so
$$
\begin{array} {rcl}
P_r(\cN_n(r;b)) &=& \frac{(n-b+r-1)!(n-1)!}{(n-b-1)!(n+r-1)!} \\
&=& \left( \frac{n-1}{n+r-1}\right) \cdots \left(\frac{n-b}{n-b+r}\right).
\end{array}
$$
This gives the bounds
$$
\begin{array} {rcl}
P_r(\cN_n(r;b)) &\ge& \left(\frac{n-b}{n-b+r}\right)^b \\
&=& \left( 1 + \frac{r}{n-b} \right)^{-b} \ge \exp(-rb/(n-b))\ ; \\
&&\\
P_r(\cN_n(r;b)) &\le& \left(\frac{n-1}{n+r-1}\right)^b \\
&=& \left( 1 - \frac{r}{n+r-1} \right)^b \le \exp(-rb/(n+r-1))\ .   
\end{array}
$$
 
\bigskip

Recall that $t=t(n)$ is defined as the least integer such that 
$P_t(\cF_n(t)) \ge 1/2$.  
Let $\w=\w(n)$ be any function tending to infinity with $n$ such that 
$t(n)/\w(n)$ takes integer values.  We shall show that 
$P_{t/\w}(\cF_n(t/\w))\to 0$ as $n\to \infty$, or equivalently that 
$P_{t/\w}(\overline{\cF_n(t/\w)}) \to 1$.  
We may assume without loss of generality that $n\ge 3$ and $\w \ge 30$, 
and so $t \ge 30$. 

We fix $n$ for the moment and consider two cases.  

\medskip

\noindent (1) \quad 
Suppose that $t = t(n) \ge 2n-1$.  

In this case, we set $b = \left \lfloor \frac{t}{2n-1}\right \rfloor$.  Note 
that our assumption on $t$ ensures that $b\ge 1$, and that our choice of $b$ 
ensures that $b(n-1)/(t-b) \le 1/2$.  Now we have
$$
\begin{array}{rcl}
P_{t-1}(\cM_n(t-1;b)) &\le& 1 - \exp(-b(n-1)/(t-b)) \\
&\le& 1- e^{-1/2} \\
&<& 1/2 \\ 
&<& P_{t-1}(\overline{\cF_n(t-1)})\ .
\end{array}
$$
Then, since $\overline{\cF_n}$ is decreasing and $\cM_n(t-1;b)$ is an initial
segment of the colexicographic order, the Clements-Lindstr\"om Theorem implies 
that
$$
\begin{array} {rcl}
P_{t/\w}(\overline{\cF_n(t/\w)}) &\ge& 
P_{t/\w}(\partial^{t-1-t/\w}\overline{\cF_n(t-1)}) \\
&\ge& P_{t/\w}(\partial^{t-1-t/\w}\cM_n(t-1;b)) \\
&=& P_{t/\w}(\cM_n(t/\w;b)) \\
&\ge& 1 - \left( \frac{t/\w}{n+t/\w-1} \right)^b.   
\end{array}
$$
If $n-1 \ge \frac{t}{\sqrt \w}$, then this gives
$P_{t/\w}(\overline{\cF_n(t/\w)}) \ge 1 - \frac{1}{\sqrt \w + 1}$.  
On the other hand, if $n-1\le \frac{t}{\sqrt \w} \le \frac{t}{5}$, then 
$\frac{t}{2n-1} \ge \frac{5n-5}{2n-1} \ge 2$, so 
$b \ge \frac{1}{2}\frac{t}{2n-1} \ge \frac{t}{5(n-1)} \ge \sqrt \w/5$, so 
$$
\begin{array}{rcl}
P_{t/\w}(\overline{\cF_n(t/\w)}) &\ge& 
1 - \left( \frac{t/\w}{t/5b+t/\w} \right)^b \\
&=& 1 - \left( 1 - \frac{\w}{\w+5b}\right)^b \\
&\ge& 1 - \exp(-b\w/(\w+5b)) \\
&\ge& 1 - \exp(\sqrt \w / 10)\ .
\end{array} 
$$

\medskip

\noindent (2) \quad
Now suppose that $t=t(n) \le 2n-2$. 

This time we set $b=\left\lceil \frac{n+t-2}{t-1} \right\rceil$, and use 
exactly the same method as in (1), but this time comparing with $\cN_n(t-1;b)$.
Indeed
$$
\begin{array}{rcl}
P_{t-1}(\cN_n(t-1;b)) &\le& \exp(-(t-1)b/(n+t-2)) \\
&\le& e^{-1} \\
&<& 1/2 \\ 
&<& P_{t-1} (\overline{\cF_n(t-1)})\ .
\end{array} 
$$
As before we deduce that 
$$
P_{t/\w}(\overline{\cF_n(t/\w)}) \ge P_{t/\w}(\cN_n(t/\w;b)) \ge 
\exp(-tb/\w(n-b))\ . 
$$
Now observe that $b \le 5n/t$, and so 
$P_{t/\w}(\overline{\cF_n(t/\w)}) \ge \exp(-6/\w)$.

\medskip

Summarizing, we see that 
$$
P_{t/\w}(\overline{\cF_n(t/\w)}) \ge \min\left\{1 - \frac{1}{\sqrt {\w(n)}+1}, 
1 - \exp(-\sqrt {\w(n)}/10),\exp(-6/\w(n)) \right\},
$$  
so $P_{t/\w}(\overline{\cF_n(t/\w)}) \to 1$ as $n\to \infty$.  

\bigskip

Now let $\w(n) \to \infty$ be such that $t(n)\w(n)$ takes integer values; we 
claim that $P_{t\w}(\overline{\cF_n(t\w)}) \to 0$ as $n\to \infty$.  We may 
assume that $n\ge 2$ and $\w(n)\ge 36$. 

As before, we fix $n$ for the moment and consider two cases. 

\medskip

\noindent (1) \quad
Suppose $t = t(n) \ge n/2$. In this case, we set 
$b = \left\lceil \frac{n+t-1}{n-1}\right\rceil$, and note that 
$b = 1 + \lceil t/(n-1) \rceil \le \frac{2t}{n-1} + \frac{2t}{n-1} = 4t/(n-1)$
and $b-1\le t\w/2$.  
Observe that
\begin{eqnarray*}
P_t(\cM_n(t;b)) & \ge & 1 - \exp(-b(n-1)/(n+t-1)) \\
	& \ge & 1 - e^{-1} \\
	& > & 1/2 \\
	& \ge & P_t(\overline{\cF_n(t)}).
\end{eqnarray*}
This implies that
$P_{t\w}(\overline{\cF_n(t\w)}) \le P_{t\w}(\cM_n(t\w;b))$, since otherwise we 
would have
$$
\left|\overline{\cF_n(t)}\right| \ge 
\left|\partial^{t\w-t}\overline{\cF_n(t\w)}\right| >
\left|\partial^{t\w-t}\cM_n(t\w;b)\right| = 
\left|\cM_n(t;b)\right|.
$$
Therefore we see that
$$
\begin{array}{rcl}
P_{t\w}(\overline{\cF_n(t\w)}) &\le& 1 - \exp(-b(n-1)/(t\w-b+1)) \\
&\le& 1 - \exp(-4t/(t\w/2)) \\
&=& 1 - \exp(-8/\w)\ .
\end{array}
$$

\medskip

\noindent (2) \quad
Now suppose $t \le (n-1)/2$.  Here we set 
$b = \left\lfloor \frac{n}{2t+1} \right\rfloor \ge 1$, and observe that
$tb/(n-b) \le 1/2$.  
We see that
\begin{eqnarray*}
P_t(\cN_n(t;b)) & \ge & \exp(-tb/(n-b)) \\
	& \ge & e^{-1/2} \\
	& > & 1/2 \\
	& \ge & P_t(\overline{\cF_n(t)})\ , 
\end{eqnarray*}
which, as before, implies that
$$
P_{t\w}(\overline{\cF_n(t\w)}) \le P_{t\w}(\cN_n(t\w;b)) \le 
\left(\frac{n-1}{n+t\w-1}\right)^b.
$$

If $n-1 \le t\sqrt \w$, this yields 
$P_{t\w}(\overline{\cF_n(t\w)}) \le 1/(1+\sqrt \w)$.  On the other hand, if
$n-1 \ge t\sqrt \w$, then 
$\frac{n}{2t+1} \ge \frac{n-1}{3t} \ge 2$, so 
$b \ge \frac{n-1}{6t} \ge \frac{\sqrt \w}{6}$, and 
$$
\begin{array} {rcl}
P_{t\w}(\overline{\cF_n(t\w)}) &\le& \left( \frac{6tb}{6tb + t\w} \right)^b \\
&=& \left( 1 - \frac{\w}{6b+\w} \right)^b \\
&\le& \exp(-\w b/(6b+\w)) \\
&\le& \exp(-\sqrt \w / 12)\ .
\end{array}
$$

Combining all the cases gives
$$
P_{t\w}(\overline{\cF_n(t\w)}) \le \max \left\{
1- \exp(8/\w(n)), \frac{1}{1+\sqrt{\w(n)}}, \exp( -\sqrt {\w(n)}/12) 
\right\},
$$
so $P_{t\w}(\overline{\cF_n(t\w)}) \to 0$ as $n\to \infty$.  

This completes the proof. 
\hfill$\pf$
\medskip

%######################################################################
%
%       SECTION:
%
\section{The Threshold for Paths}\label{PathThresh}
For convenience we will assume that all logarithms and roots
take on the value of their nearest integer.
Let $P_n$ be the path on $n$ vertices, $c$ be any constant less than $1/\sq$,
$u=c\sln$, $t=n2^u$, and $p=(1+\e)2^u \ln n$, for some $\e>0$. 
In addition, let $D_i$ be the random variable
which is the number of pebbles on vertex $i$ when a random
distribution of $t$ pebbles is selected, and let $S_i$ denote the
event that $D_i\le p$ and $T_i=\neg S_i$. 
Finally, we will also need the following random variables
\begin{equation}\label{y+}
Y_i^+=\sum_{l=i}^n{D_l\over 2^{l-i}}
\end{equation}
and 
\begin{equation}\label{y-}
Y_i^-=\sum_{l=1}^i{D_l\over 2^{i-l}}\ .
\end{equation}

\begin{lemma}\label{union}
Let $p,t,n,u$ be as above. 
Then $\Pr [T_1\cup T_2\cup\cdots\cup T_n]\rar 0$.
\end{lemma}

\noindent
{\it Proof.}
For $1\le i\le n$, we have
$$\Pr [T_i]
	={{n\pebb t-p}\over {n\pebb t}}
	={(t-p+1)\cdots (t)\over (t+n-p)\cdots (t+n-1)}$$
$$\le\left({t\over t+n}\right)^p
	\le e^{-np/(t+n)}\ .$$
Therefore,
$$\Pr [T_1\cup\cdots\cup T_n]\le ne^{-p/(2^u+1)}=
\exp\left( \ln n \left( 1- {2^u \over 2^u+1}(1+\e) \right) \right) \rar 0.$$
\hfill$\pf$
\medskip

Note that the role of $u$ is not important in the above lemma.
However, we stated the lemma as such because of its use in the proof of
Theorem \ref{Thpaths}.

For the next lemma, we will need $m=\stln$, and $k=\lf n/m\rf$.
Partition $P_n$ into consecutively disjoint paths (blocks)
$B_1,\ldots,B_k$ of lengths $m$ or $m+1$.
For a randomly chosen distribution of $t$ pebbles on $P_n$ we denote by
$E_i$ the event that block $B_i$ contains no pebbles (is empty), and set
$F_i=\neg E_i$.

\begin{lemma}\label{mainlower}
Let $m,k$ be as above and let $s=s(n)$ be such that $n\gg ms^me^{m/s}$.
If $t=ns$ then $\Pr [F_1\cap\cdots\cap F_k]\rar 0$.
\end{lemma}

\noindent
{\it Proof.}
We will apply the second moment method. 
Let $X_i$ be the indicator variable of $E_i$ and let $X=\sum_{i=1}^k X_i$. 
We have
$$\Pr(X=0)\le{\sigma^2(X)\over (\bE(X))^2}$$
and
\begin{eqnarray*}
\sigma^2(X)&=&\bE (X^2)-(\bE(X))^2\\
&&\\
&=&\sum_{i,j}\bE (X_iX_j)-\sum_{i,j}\bE (X_i)\bE (X_j)\\
&&\\
&\le&\sum_i\bE (X_i^2)\ ,
\end{eqnarray*}
where the last inequality follows from the fact that 
$\bE (X_iX_j)\le\bE (X_i)\bE (X_j)$ for $i\ne j$.
Therefore,
$$\sigma^2(X)\le\sum_i\bE (X_i^2)=\sum_i\bE (X_i)=\bE (X)\ .$$
We also have 
\begin{eqnarray*}
\bE(X)&=&\left({n\over m}\right){{t+n-m-1\choose t}\over {t+n-1\choose t}}\\
&&\\
&\ge&\left({n\over m}\right)\left({n-m\over t+n-m}\right)^m\\
&&\\
&\sim&\left({n\over ms^m}\right)e^{-m^2/n -m/s}\\
&&\\
&\rar&\inf\ .
\end{eqnarray*}
Thus 
$\Pr[X=0]\rar 0.$
\hfill$\pf$

\medskip
\noindent
{\it Proof of Theorem \ref{Thpaths} (a).}
We first remark that it is easy to check that $s=2^u$ satisfies
hypothesis of Lemma \ref{mainlower}.
Next, for each block $B=[i,j]$, let $Y_B^-=Y_i^-$ and $Y_B^+=Y_j^+$. 
Then Lemma \ref{union} and Lemma \ref{mainlower} imply that with probability 
tending to one there is a block $B$ of length $m$ which is empty and such 
that both $Y_B^-$ and $Y_B^+$ are less than or equal to $(1+\e)2^u \ln n $.
Since $(1+\e)2^u \ln n <2^{m/2}$ for large enough $n$, 
there is a vertex of $B$, namely the center vertex, 
to which it is not possible to pebble using $Y_B^-$ and $Y_B^+$. 
\hfill$\pf$
\medskip

We also shall prove that $th(P_n) \subseteq O(n2^{2\sln})$. 
The argument is a modification of an idea from \cite{CzyEHK}.
We will use the fact \cite{CzyEHK} that any distribution of at least $2^m$
pebbles on a connected graph of $m+1$ vertices is solvable.

\medskip
\noindent
{\it Proof of Theorem \ref{Thpaths} (b).}
Let $m=\sln$, $k=\lf n/m\rf$, $s=2^{2m+2}$, $t=ns$, and 
partition $P_n$ into consecutively disjoint paths (blocks)
$B_1,\ldots,B_k$, each with $m$ or $m+1$ vertices.
Let $Z_i$ denote the event that the $i$th block contains less than 
$2^m$ pebbles. 
We will show that $\Pr [Z_1\cup\cdots\cup Z_k]\rar 0$, which, by the remark
preceding the proof, implies the result.  

We have that
$$\Pr [Z_1\cup\cdots\cup Z_k] 
	\le k\sum_{i=0}^{2^m-1}
	{{i+m-1\choose i}{t+n-m-i-1\choose t-i}\over {t+n-1\choose t}}\ .$$
Since 
${t+n-m-i-1\choose t-i} 
	\le\left({t\over t+n-m-1}\right)^i{t+n-m-1\choose t}
	\le {t+n-m-1\choose t} 
	\le\left({n\over n+t}\right)^m {t+n-1\choose t}$, we have
\begin{equation}\label{eq5}
{t+n-m-i-1\choose t-i} 
	\le\left({n\over n+t}\right)^m{t+n-1\choose t}\ .
\end{equation}
For $0\le i < 2^m$, we have 
\begin{equation}\label{eq6}
{i+m-1\choose i}\le 2^m2^{m^2}\ .
\end{equation}
Using (\ref{eq5}) and (\ref{eq6}), we see that
$$\Pr [Z_1\cup\cdots\cup Z_k] 
	\le k2^{2m+m^2}\left({n\over n+t}\right)^m 
	\le {n\over m}\left({1\over s}\right)^m 2^{2m+m^2}
	={1\over m}\rar 0\ .$$
\hfill$\pf$
\medskip

%######################################################################
%######################################################################
%
%       BIBLIOGRAPHY:
%
\bibliographystyle{plain}
%   

%######################################################################
%######################################################################
%
%       END OF PAPER
%
\end{document}